\documentclass[letterpaper,11pt]{article}
\usepackage{amsmath, amsthm, amssymb}    
\usepackage{bridges}			
\usepackage{graphicx}			
\usepackage[colorlinks=true, urlcolor=blue, citecolor=black, linkcolor=black]{hyperref}  
\usepackage{subcaption}			
\usepackage{xcolor}

\usepackage{booktabs}
\usepackage{mathtools}
\usepackage{wrapfig}

\usepackage{pgfplots}

\usepackage{tikz}
\usetikzlibrary{perspective,plotmarks}

\definecolor{mvj_blue}{HTML}{333676}
\definecolor{mvj_pink}{HTML}{611047}
\definecolor{mvj_orange}{HTML}{804c15}
\definecolor{mvj_gold}{HTML}{aaa239}
\definecolor{mvj_green}{HTML}{116416}

\tikzset{every picture/.style={ultra thick}}

\newcommand{\tetrahedron}[1]{
\begin{tikzpicture}[3d view={-45}{150}]
    \coordinate (v0) at (0,0,0);
    \coordinate (v1) at (1,1,0);
    \coordinate (v2) at (0,1,1);
    \coordinate (v3) at (1,0,1);
    \foreach \x/\y in {#1} {
     \draw [mvj_blue] (v\x) -- (v\y);
    }
\end{tikzpicture}
}

\newcommand{\cube}[1]{
\begin{tikzpicture}[3d view={-45}{150}]
    \coordinate (v0) at (0,1,0);
    \coordinate (v1) at (0,1,1);
    \coordinate (v2) at (0,0,1);
    \coordinate (v3) at (0,0,0);
    \coordinate (v4) at (1,1,0);
    \coordinate (v5) at (1,1,1);
    \coordinate (v6) at (1,0,1);
    \coordinate (v7) at (1,0,0);
    \foreach \x/\y in {#1} {
     \draw [mvj_pink] (v\x) -- (v\y);
    }
\end{tikzpicture}
}

\newcommand{\octahedron}[1]{
\begin{tikzpicture}[3d view={-45}{150}]
    \coordinate (v0) at (0,0,1);
    \coordinate (v1) at (1,0,0);
    \coordinate (v2) at (0,1,0);
    \coordinate (v3) at (-1,0,0);
    \coordinate (v4) at (0,-1,0);
    \coordinate (v5) at (0,0,-1);
    \foreach \x/\y in {#1} {
     \draw [mvj_orange] (v\x) -- (v\y);
    }
\end{tikzpicture}
}

\newcommand{\dodecahedron}[1]{
\begin{tikzpicture}[3d view={-45}{150}]
\coordinate (v0) at (1.00,1.00,1.00);
\coordinate (v1) at (0.00,1.62,0.62);
\coordinate (v2) at (-1.00,1.00,1.00);
\coordinate (v3) at (-0.62,0.00,1.62);
\coordinate (v4) at (0.62,0.00,1.62);
\coordinate (v5) at (1.62,0.62,0.00);
\coordinate (v6) at (0.00,1.62,-0.62);
\coordinate (v7) at (-1.62,0.62,0.00);
\coordinate (v8) at (-1.00,-1.00,1.00);
\coordinate (v9) at (1.00,-1.00,1.00);
\coordinate (v10) at (1.00,1.00,-1.00);
\coordinate (v11) at (-1.00,1.00,-1.00);
\coordinate (v12) at (-1.62,-0.62,0.00);
\coordinate (v13) at (0.00,-1.62,0.62);
\coordinate (v14) at (1.62,-0.62,0.00);
\coordinate (v15) at (1.00,-1.00,-1.00);
\coordinate (v16) at (0.62,0.00,-1.62);
\coordinate (v17) at (-0.62,0.00,-1.62);
\coordinate (v18) at (-1.00,-1.00,-1.00);
\coordinate (v19) at (0.00,-1.62,-0.62);
    \foreach \x/\y in {#1} {
     \draw [mvj_gold] (v\x) -- (v\y);
    }
\end{tikzpicture}
}

\newcommand{\icosahedron}[1]{
\begin{tikzpicture}[3d view={-45}{150}]
\coordinate (v0) at (0.00,1.00,1.62);
\coordinate (v1) at (0.00,-1.00,1.62);
\coordinate (v2) at (1.62,0.00,1.00);
\coordinate (v3) at (1.00,1.62,0.00);
\coordinate (v4) at (-1.00,1.62,0.00);
\coordinate (v5) at (-1.62,0.00,1.00);
\coordinate (v6) at (-1.00,-1.62,0.00);
\coordinate (v7) at (1.00,-1.62,0.00);
\coordinate (v8) at (1.62,0.00,-1.00);
\coordinate (v9) at (0.00,1.00,-1.62);
\coordinate (v10) at (-1.62,0.00,-1.00);
\coordinate (v11) at (0.00,-1.00,-1.62);
    \foreach \x/\y in {#1} {
     \draw [mvj_green] (v\x) -- (v\y);
    }
\end{tikzpicture}
}

\title{Incomplete Open Platonic Solids}
\author{Mikael Vejdemo-Johansson}
\date{[Draft as of \today]}

\begin{document}

\maketitle

\thispagestyle{empty}

\begin{abstract}
\hfill
\tetrahedron{0/1,0/2,0/3,1/2}
\hfill
\cube{0/3,0/4,1/2,1/5,2/6,3/7,4/5}
\hfill
\octahedron{0/1,0/2,0/3,0/4,1/2,1/4,1/5,2/3,3/4}
\hfill
\dodecahedron{0/1,0/4,0/5,1/2,1/6,2/3,2/7,3/4,3/8,4/9,5/10,5/14,6/10,7/11,7/12,8/12,9/13,9/14,10/16,11/17,12/18,13/19,15/16,15/19,17/18}
\hfill
\icosahedron{0/1,0/2,0/3,0/4,0/5,1/2,1/6,1/7,2/3,2/7,2/8,3/4,3/8,3/9,4/5,4/9,4/10,5/6,5/10,6/10,6/11,7/8,7/11,8/9,8/11}
\hfill

    Sol LeWitt famously enumerated all the incomplete open cubes, finding 122 of these connected, non-planar subsets of the edges of the cube.
    Since then, while several projects have revisited the cube enumeration, no such enumeration has been published for any other interesting solid.

    In this paper we present work on enumerating all the incomplete open platonic solids, finding 6 tetrahedra, 122 cubes (just like LeWitt), 185 octahedra, 2\,423\,206 dodecahedra and 16\,096\,166 icosahedra.
\end{abstract}

\section{Introduction}

\begin{wrapfigure}{r}{6.5cm}
\includegraphics[width=\linewidth]{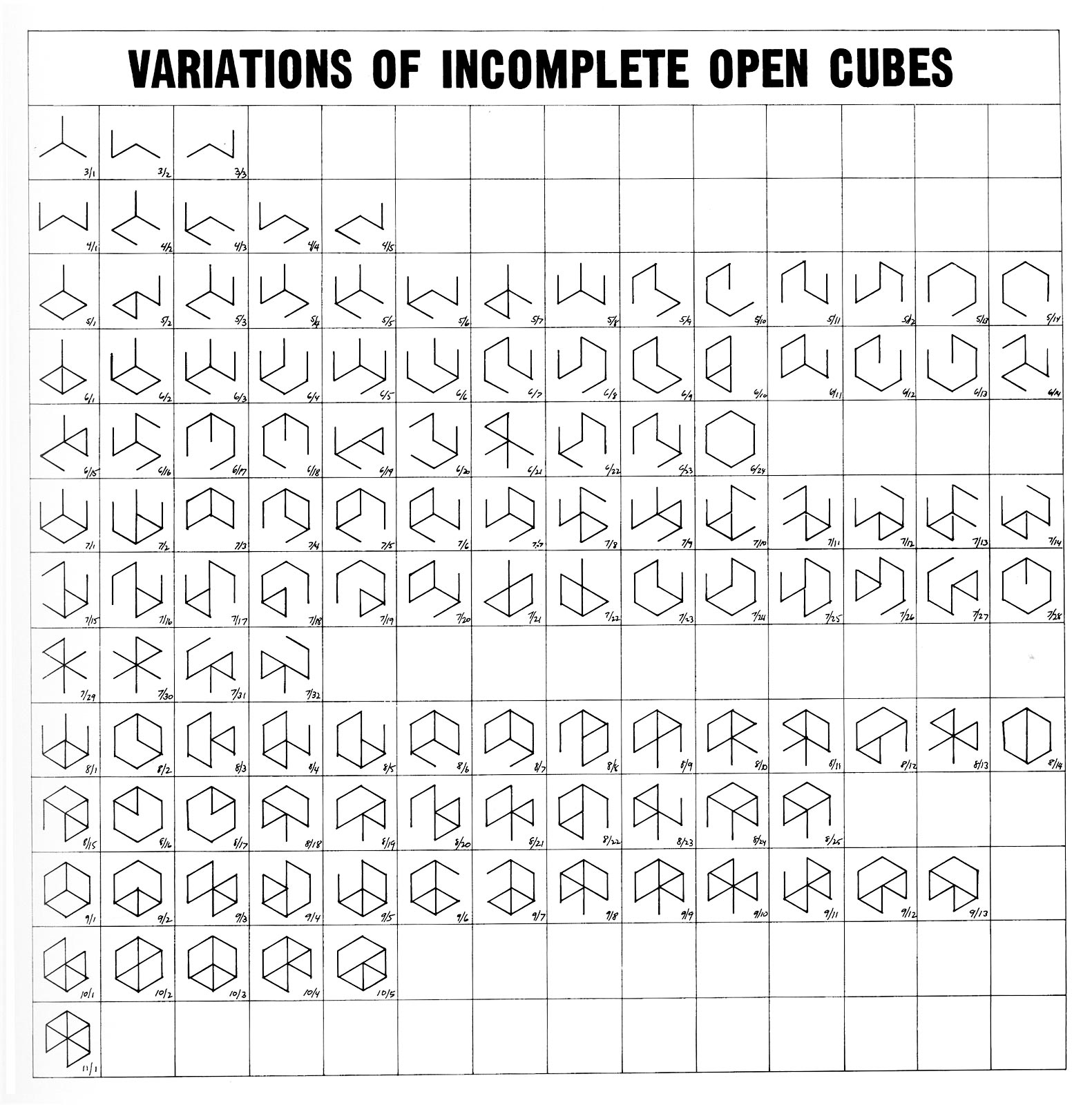}
\caption{Exhibition poster for Sol LeWitt's Variations on Incomplete Open Cubes from a show at the John Weber Gallery in 1974.}\label{fig:show-poster}
\end{wrapfigure}

Sol LeWitt's \emph{Variations on Incomplete Open Cubes} has drawn several waves of attention from mathematicians: the task that the artist set himself is deeply mathematical in nature, and can be well understood using symmetry groups and group actions. 
Nevertheless, the artist himself worked out the full list by hand and arrived at a correct count (later verified unpublished by Dr. Erna Herney and Prof. Arthur Babakhanian, and published with the location of one mistake -- including a rotation instead of a reflection -- by Rozhovskaya and Reb \cite{rozhkovskaya2015list} and \cite{reb2013analysis}).
The project has also been the inspiration for Rob Weychert's visualization \emph{Incomplete Open Cubes Revisited} \cite{Weychert_2018}, noting that excluding the empty cube and the full cube, the remaining cubes total a count of 4094 if one includes both planar, disconnected and rotated into the count.
More recently, Alison Pitt \cite{Pitt_2023} created an interactive Tableau visualization, and 3blue1brown released a video going in depth on the connections to group theory and group actions and using Burnside's Lemma to validate Sol LeWitt's count \cite{3b1b}.

As far as I have been able to tell, there has been no work going into analogous constructions on other solids.
From LeWitt's side, it is understandable -- he asserted \textit{``The cube is the best form to use as a basic unit for any more elaborate function, the grammatical device from which the work may proceed.''} in 1966 (as described by the Paula Cooper Gallery \cite{coopergallery2013}), and with that approach would not have been interested in other shapes.
Here, we seek to fill the gap by computing the complete enumerations for all five platonic solids, with a methodology that is easy to generalize to other symmetric wireframes.

\section{Ruleset: what is an incomplete open polyhedron?}

Sol LeWitt defines an incomplete open cube as a \emph{connected} and \emph{non-planar} subset of the edges of the cube and enumerates them up to rotational equivalence: mirror images are allowed, but the same shape in a different rotation is considered not to be a new shape.

Following that we will work in this piece with the notion of an incomplete open polyhedron to be a \emph{connected} and \emph{non-planar} \emph{proper} and \emph{non-empty} subset of the set of edges of the polyhedron.

\section{Enumeration}

To enumerate all the incomplete open polyhedra, I ended up using a process in several steps, drawing on the data from the SymPy python package to extract information on the polyhedra, and then translating the vertex permutations from SymPy into edge permutations and using bit manipulation in C++ to then find single representatives for each orbit of the symmetry group action.

The code used for this project is made available at \url{https://github.com/michiexile/incomplete-open}\footnote{Project code is at: DOI:10.5281/zenodo.18749808}.

\subsection{Subset encoding}

Since all the polyhedra under consideration have fewer than 32 edges\footnote{Tetrahedron: 6, cube and octahedron: 12, dodecahedron and icosahedron: 30}, we were able to use bit patterns in a standard long integer (using \texttt{uint32\_t} to be specific) to represent subsets of edges for all the polyhedra.
The representation in SymPy has each polyhedron represented by integer-labeled vertices, and with edge sets corresponding to these vertices. We represent these edge sets by sorting each edge, and then sorting the entire set lexicographically: that way each edge has a unique list index position that we can use to adress each edge set by the bits in a 32-bit word.

\subsection{Symmetry groups and their actions on edges and subsets}

SymPy helpfully encodes the vertex action of each symmetry group: each group element $g$ moves each vertex $v$ to some vertex $gv$.
Each edge is determined as a pair of vertices: $[v_1,v_2]\xmapsto{g}[gv_1,gv_2]$; by sorting the result we get a new edge in the same format as in our original listings -- so we can lookup the index of the edge that each edge maps to.

This let's us translate vertex actions into edge actions, for example, in the octahedral group acting on the cube, the vertex permutation $(0\ 6)(1\ 7)(2\ 3)(4\ 5)$ yields the edge permutation $(0\ 11)(1\ 6)(2\ 10)(3\ 7)(4\ 9)$.

The set of all possible images $\{gs : g\in G, s\in S\}\subseteq S$ for a group $G$ acting on a set $S$ is called an \emph{orbit}, and finding edgesets up to rotation corresponds exactly to finding representatives of each orbit of the rotation group action on the set of edge sets.
As soon as we have a well-ordering on the set $S$, we can easily define a notion of canonical representatives to be the minimal elements of each orbit.
We will use this to create a process that quickly rejects edgesets that cannot be minimal elements of their orbits, leaving only the minimal representatives at the end of the process.

\subsection{Putting everything together computationally}

For the smaller cases: tetrahedron, cube and octahedron, we could enumerate all possible edgesubsets and compute with each directly -- it takes a little time, but is quite feasible since $2^{12}=4\,096$. $2^{30}$ however is much larger: for dodecahedra and icosahedra we would need to work through $1\,073\,741\,824$ subsets in total. The count of permutations to consider explodes similarly fast.

However, enumerating all possible subsets of 8 elements is considerably easier: there are 256 of them, and a quick computation can establish a lookup table for where each edge should go. With one of these tables we can take care of the tetrahedron, two tables (one for the first 8 edges, a second one for the 4 remaining edges) can handle the cube and octahedron, and four tables comfortably handles the 30 edges for the dodecahedron and icosahedron.

Each table takes an index number $0\leq i<256$ corresponding to a particular subset of the corresponding run of edges. 
At that index, the table then stores a 32-bit word with the destination bits for all the bits set in the index number. 
Since each permutation is a bijection, these 32-bit words will be disjoint, and we can use a bit-wise OR operation to combine the permutation action on each run of 8 edges into a complete description of the edge permutation.

With these lookup tables in place, the search for a canonical representative for each orbit is fast: a given subset is the canonical representative if all permutations of its bits are larger (as 32-bit unsigned integers) than the representative itself. 
While we are testing, we can stop early the moment we find something smaller than the element we are testing.
Since testing of each subset is independent of each other subset, this task is what is known as an \emph{embarrassingly parallel} problem: it is very easy to parallelize.

Once we have orbit-minimal representatives, we can then remove the ones that do not fit the rule sets. 
There are two steps we need to handle: connectivity, and non-planarity.

\paragraph{Connected} There is a component of the Boost library for C++ programming that computes graph properties -- we use this to count connected components of each edge set under consideration.

\paragraph{Planar} With some trial and error, we found vertex coordinates that are compatible with the edge structures encoded in SymPy.
With these vertex coordinates, we can test for planarity using a normal vector: use the first three vertices $p_0,p_1,p_2$ to construct vectors $v_1=p_1-p_0$ and $v_2=p_2-p_0$. The normal to the plane spanned by $v_1,v_2$ is $v_1\times v_2$, and for another point $p$, the dot product $(v_1\times v_2)\cdot(p-p_0)$ is zero exactly when $p$ falls in the plane spanned by $p_0,p_1,p_2$.

After exporting permutations, edge lists and vertex positions from SymPy and our reconstruction, the full enumeration took the following times\footnote{Reporting the \texttt{real} time from using the \texttt{time} utility on Ubuntu Linux} on a compute server with 32 cores, all fully dedicated to the numeration: tetrahedron (0.056), cube (0.079 seconds), octahedron (0.084 seconds), dodecahedron (107 minutes and 4.362 seconds) and icosahedron (109 minutes and 5.067 seconds).

\subsection{All the incomplete open platonic solids}

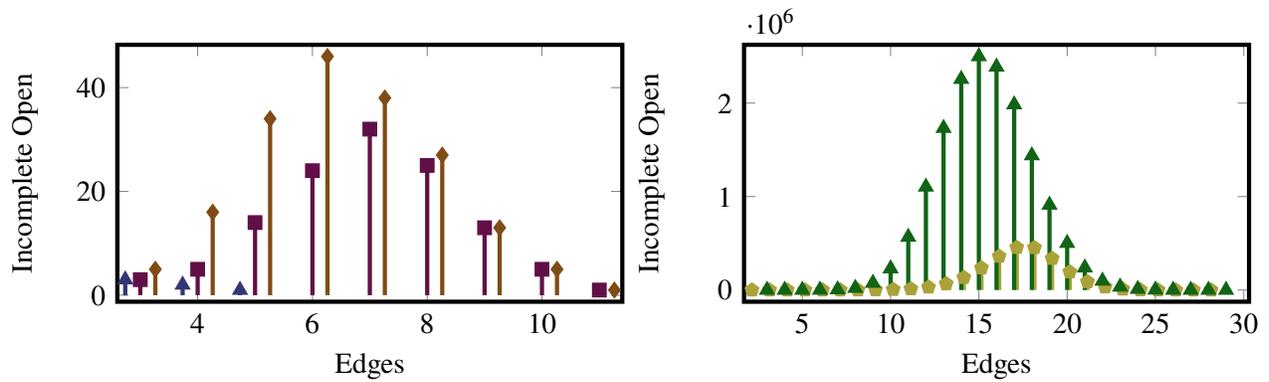
\begin{figure}
    \begin{tikzpicture}
    \begin{axis}[
    clip bounding box=upper bound,
    enlargelimits=0.05,
    ycomb,
    ylabel=Incomplete Open,
    xlabel=Edges,
    height=5cm, width=0.5\linewidth
    ]
    \addplot+[color=mvj_blue,mark options={mark=triangle*},scale=2.0,xshift=-0.2cm] coordinates 
    {(3,3) (4,2) (5,1)};
    \addplot+[color=mvj_pink,mark options={mark=cube*}] coordinates
    {(3,3) (4,5) (5,14) (6,24) (7,32) (8,25) (9,13) (10,5) (11,1)};
    \addplot+[color=mvj_orange,mark options={mark=diamond*},xshift=+0.2cm] coordinates {(3,5) (4,16) (5,34) (6,46) (7,38) (8,27) (9,13) (10,5) (11,1)};
    \end{axis}
\end{tikzpicture}
\begin{tikzpicture}
    \begin{axis}[
    clip bounding box=upper bound,
    enlargelimits=0.05,
    ycomb,
    ylabel=Incomplete Open,
    xlabel=Edges,
    height=5cm, width=0.5\linewidth
    ]
    \addplot+[color=mvj_gold,mark options={mark=pentagon*},scale=2.0,xshift=-0.2cm] coordinates 
    {(3,3) (4,5) (5,18) (6,43) (7,119) (8,300) (9,818) (10,2083) (11,5357) (12,13078) (13,30674) (14,66723) (15,133347) (16,236182) (17,360834) (18,455307) (19,452799) (20,338011) (21,193929) (22,88217) (23,32545) (24,9834) (25,2408) (26,482) (27,78) (28,11) (29,1)};
    \addplot+[color=mvj_green,mark options={mark=triangle*}] coordinates
    {(3,10) (4,45) (5,234) (6,1080) (7,4936) (8,20408) (9,74825) (10,229386) (11,567132) (12,1103778) (13,1731869) (14,2256123) (15,2498922) (16,2387598) (17,1984459) (18,1439268) (19,910702) (20,501388) (21,238981) (22,97889) (23,34112) (24,10025) (25,2423) (26,483) (27,78) (28,11) (29,1)};
    \end{axis}
\end{tikzpicture}
\vspace{-1.75cm}
    \caption{Number of different incomplete polyhedra by edge-count for each. The left graph shows the tetrahedron ({\protect\tikz{\protect\node[mvj_blue]{\protect\pgfuseplotmark{triangle*}};}}), the cube ({\protect\tikz{\protect\node[mvj_pink]{\protect\pgfuseplotmark{square*}};}}) and the octahedron ({\protect\tikz{\protect\node[mvj_orange]{\protect\pgfuseplotmark{diamond*}};}}) and the right graph shows the dodecahedron ({\protect\tikz{\protect\node[mvj_gold]{\protect\pgfuseplotmark{pentagon*}};}}) and the icosahedron ({\protect\tikz{\protect\node[mvj_green]{\protect\pgfuseplotmark{triangle}};}}).}\label{fig:counts}
\end{figure}

After all of this, we arrive at a complete enumeration of the incomplete open platonic solids.
In \autoref{tab:enumeration}, we see the resulting counts from the enumeration. To browse the complete collection of solids, we invite you to the viewer at \url{https://incomplete-open.mikael.johanssons.org}.

\begin{table}
    \centering
    \begin{tabular}{llr|llr}
    \toprule
        Solid & Symmetry group & Total count &
        Solid & Symmetry group & Total count \\ \midrule
        Tetrahedron & Tetrahedral: $A_4$ & 6 \\
        Cube & Octahedral: $S_4$ & 122 &
        Octahedron & Octahedral: $S_4$ & 185 \\
        Dodecahedron & Icosahedral: $A_5$ & 2\,423\,206 &
        Icosahedron & Icosahedral: $A_5$ & 16\,096\,166 \\
        \bottomrule
    \end{tabular}
    \caption{Summary of the complete enumeration of incomplete open platonic solids.}
    \label{tab:enumeration}
\end{table}

\section{Discussion}

Looking at the finished enumeration, there are many questions that emerge -- some at first glance, others when looking deeper. One of the more fascinating ones is the observation that we get different counts -- sometimes spectacularly different counts -- from dual polyhedra: the cube and octahedron share edge-count and symmetry group, the dodecahedron and icosahedron also share edge-count and symmetry group, yet the octahedron enumeration is 51\% larger than the cube enumeration and the icosahedron has $6.6\times$ more incomplete open polyhedra than the dodecahedron.

The key to these differences is in connectivity (and in node degrees): a connected set of edges in the cube might become disconnected when dualized\footnote{the effect of dualizing a polyhedron is to rotate each edge by 90\textdegree.}, which accounts for most of the difference between cubes and octahedra. For the dodecahedra and icosahedra, I believe the difference is largely due to that each vertex is incident to $3$ edges on a dodecahedron, but $5$ edges on an icosahedron: there are many more paths on the icosahedral edges than on the dodecahedral ones, which reflects in the combinatorial explosion of connected edge sets.

\bibliographystyle{bridges}
\bibliography{lewitt}

@article{rozhkovskaya2015list,
  title={Is the list of incomplete open cubes complete?},
  author={Rozhkovskaya, Natasha and Reb, Michael},
  journal={Nexus Network Journal},
  volume={17},
  number={3},
  pages={913--925},
  year={2015},
  publisher={Springer}
}

@mastersthesis{reb2013analysis,
  title={Analysis of variations of Incomplete open cubes by Sol Lewitt},
  author={Reb, Michael Allan},
  year={2013},
  type={Bachelor's Thesis},
  school={Kansas State University}
}

@misc{Weychert_2018, 
  title={Incomplete open cubes revisited}, 
  url={https://cubes-revisited.art/}, 
  journal={Incomplete Open Cubes Revisited}, 
  author={Weychert, Rob}, 
  year={2018}
}

@misc{ Pitt_2023, 
  title={An exploration of incomplete open cubes}, 
  url={https://www.alisonpitt.com/blog/an-exploration-of-incomplete-open-cubes}, 
  journal={Alison Pitt}, 
  author={ Pitt, Alison}, 
  year={2023}, 
  month={Sep}
}

@misc{3b1b, 
  title={Exploration and Epiphany}, 
  url={https://www.youtube.com/watch?v=_BrFKp-U8GI}, 
  journal={3Blue1Brown}, 
  publisher={YouTube}, 
  author={Dancstep, Paul}, 
  year={2025}, 
  month={Sep}
}

@misc{coopergallery2013, 
  title={Sol LeWitt -- Cubic Forms},
  url={https://www.paulacoopergallery.com/exhibitions/sol-lewitt13}, 
  publisher={Paula Cooper Gallery}, 
  year={2013}
}

\end{document}